\documentclass[a4paper,10pt]{article}

\usepackage{amsfonts}
\usepackage{amssymb}
\usepackage{amsthm}
\usepackage{amsmath}
\usepackage{latexsym}
\usepackage[all]{xy}
\usepackage{mathrsfs}
\usepackage{graphicx}

\theoremstyle{plain}\newtheorem{thm}{Theorem}[section]
\theoremstyle{plain}\newtheorem{lem}[thm]{Lemma}
\theoremstyle{plain}\newtheorem{prop}[thm]{Proposition}
\theoremstyle{plain}\newtheorem{cor}[thm]{Corollary}
\theoremstyle{definition}\newtheorem{defn}[thm]{Definition}
\theoremstyle{remark}\newtheorem{rem}[thm]{REMARK}
\theoremstyle{remark}

\title{Differential equations, difference equations and algebraic relations: An extension to a theorem of Compoint}
\author{Camilo Sanabria Malag\' on\footnote{partially supported by NSF grant CCF 0901175 and CCF 0952591.}}
\date{}

\begin{document}

\maketitle

\begin{abstract}
Let $C$ be an algebraically closed field and $X$ a projective curve over $C$. Consider an ordinary linear differential equation, or a linear difference equation, with coefficients in the field of rational functions of $X$, and assume that its Galois Group $G$ has finite determinant group and is reductive. In this context, the ideal of algebraic relations satisfied by a full system of solutions is generated by the  $G$-invariants it contains. This result extends a theorem of E. Compoint \cite{compoint}.
\end{abstract}

\section{Introduction}
Let $C$ be an algebraically closed field and let $X$ be a projective curve over $C$ with $K=C(X)$ the field of rational functions on $X$.

We fix an operator
\begin{eqnarray*}
v: K & \longrightarrow & K\\
\end{eqnarray*}
which could be either a \emph{derivation}, i.e. an additive map satisfying the Leibnitz rule, or an automorphism. The \emph{constants} of $K$ are the $f\in K$ such that $v(f)=0$, if $v$ is a derivation; or $v(f)=f$, if $v$ is an automorphism. We assume now that $v$ is such that the field of constants is $C$.

We take a matrix equation of the form:
\begin{eqnarray}\label{eqn}
v(y^i)=a^i_jy^j,\quad (a^i_j)_{1\le i,j\le n}\in M_{n\times n}(K).
\end{eqnarray}
If $v$ is a derivation, such equation is called a \emph{differential equation}; whereas if $v$ is an automorphism, the equation is a \emph{difference equation}. In the latter case we impose the condition that $(a^i_j)$ is invertible. A \emph{full system of solutions} of (\ref{eqn}) is an invertible $n\times n$ matrix $(f^i_j)$ such that
\[
v(f^i_j)=a^i_kf^k_j,\quad \forall i,j\in\{1,\ldots,n\}.
\]
We pose the following question:

 \emph{What are the algebraic relations among the $\{f^i_j\}_{1\le i,j\le n}$? in other words, what are the polynomials $P(X^i_j)\in K[X^i_j]_{1\le i,j\le n}$ such that $P(f^i_j)=0$?}

The question of effectively describing the algebraic relations satisfied by the solutions of (\ref{eqn}) is related to the question of determining the Galois group $G$ of the equation. In its turn, the Galois group may tell us some useful information about the solutions like whether these are algebraic, or whether these are liouvillian. In the differential setting, liouvillian means that one can craft the solutions by iterating integration, exponentiation of integrals or taking roots of polynomials; in the difference setting the meaning depends on the chosen automorphism. For example if $K=C(x)$ and $v(x)=x+1$, it means that the equation is solvable in terms of liouvillian sequences \cite{hensing}.

  Furthermore, if the solution satisfy specific kind of algebraic relations, we may be able to use a particular solving method or to get a more accurate description of the solution \cite{singer88} \cite{nguyen}.

  So, although our answer to the question above will be given using the Picard-Vessiot extension of (\ref{eqn}), and the result in this paper can be applied to describe this extension; our result is interesting because of the concequences algebraic relations have on the solutions.

  Here the Picard-Vessiot extension of (\ref{eqn}) is the minimal \emph{differential} or \emph{difference} (depending on the setting) \emph{$K$-algebra} containing a full system of solutions but with no new contant. The adjectives differential and difference mean that the respective operator extends to this K-algebra.

  Under appropriate hypothesis on the Galois group of the equation, the problem in the differential setting has been completely solved by Compoint~\cite{compoint}. For, theoretically, to obtain a full basis of $G$-invariant  polynomials one can rely on the algorithm by M. van Hoeij and J.A. Weil~\cite{We}, together with the bounds obtained in ~\cite{SC}; and Compoint's solution asserts that all algebraic relations are written in terms of $G$-invariants (Theorem \ref{thmcom} below). In his paper Compoint also asserts that his solution fails if we relax his hypothesis.

  However by modifying his argument, we are able to generalize his result to a larger class of equations (Theorem \ref{thm}). In particular we extend the result to the difference setting; where the computational hurdle encountered by our solution to the problem is the absence of faster algorithms for obtaining the group invariants (like van Hoeij-Weil's in the differential setting). So far, in the generic difference setting, one has to first obtain the symmetric power of the equation (cf. Remark \ref{cat}) and then find rational solutions, which is a rather slow process. As with the interpretation of solvability of the Galois group, the methods to find rational solutions depend on the automorphism $v$ (eg. \cite{hoeij}).

  Nevertheless, the similarities between the differential Galois theory and the difference counterpart allows us to give an unified argument to answer the question. A very concise introduction to both setting is in \cite{etingof}, a more detailed exposition of each may be found on \cite{SvdP0} and \cite{SvdP}. 

 Our proof will follow F. Beukers' argument~\cite{beukers}. We treat the problem in terms of differential and difference modules.

\section{Setting}

\begin{defn}
Assume $v$ is a derivation. A \emph{differential module} $M$ over $K$ is a $K$-vector space together with an additive map $\partial$ such that:
\[
\partial fm=v(f)m+f\partial m
\]
for every $f\in K$ and $m\in M$. If $\partial m=0$ we say that $m$ is horizontal.
\end{defn}

\begin{defn}
Assume $v$ is an automorphism. A \emph{difference module} $M$ over $K$ is a $K$-vector space together with an automorphism $\Phi$ such that:
\[
\Phi fm=v(f)\Phi m
\]
for every $f\in K$ and $m\in M$. If $\Phi m=m$ we say that $m$ is horizontal.
\end{defn}

\begin{rem}\label{rem1}
It is quite evident from the definition that, on both settings, the collection of horizontal elements forms a vector space over the field of constants. Given a differential module $M$ and a basis $e_1,\ldots,e_n$, if $\partial e_j=-a^i_je_i$, to solve the equation
\[
\partial m=0
\]
in this basis amounts to solve the system of equations:
\[
v(x^i)=a^i_jx^j.
\]
Indeed if $m=x^ie_i$, then
\begin{eqnarray*}
\partial m & = &  v(x^i)e_i- x^ja^i_je_i\\
           & = & \left( v(x^i)-a^i_jx^j \right)e_i.
\end{eqnarray*}
On the other hand, given a difference module $M$ and a basis $e_1,\ldots,e_n$, if $\Phi e_j= a^{-1}{}^i_je_i$ (where $a^i_ka^{-1}{}^k_j=\delta^i_j$), to solve the equation 
\[
\Phi m=m
\]
in this basis amounts to solve the system of equations:
\[
v(x^i)=a^i_jx^j.
\]
Indeed if $m=x^ie_i$, then
\begin{eqnarray*}
\Phi m & = &  v(x^i)a^{-1}{}^j_ie_j
\end{eqnarray*}
so $\Phi m=m$ is equivalent to
\begin{eqnarray*}
v(x^i)a^{-1}{}^j_ie_j & = & x^je_j.
\end{eqnarray*}
\end{rem}

\begin{rem}\label{cat}
Let $M_1$ and $M_2$ be two differential modules (resp. difference modules) and denote their respective operators by $\partial_1$ and $\partial_2$ (resp. $\Phi_1$ and $\Phi_2$). We turn $M_1\otimes_K M_2$ into a differential (resp. difference) module by operating as follows:
\begin{eqnarray*}
m_1\otimes m_2 & \longmapsto & \partial_1m_1\otimes m_2+m_1\otimes \partial_2m_2\\
\textrm{(resp.}\quad m_1\otimes m_2 & \longmapsto & \Phi_1m_1\otimes \Phi_2m_2\quad\textrm{)}
\end{eqnarray*}
In this manner, starting with a differential (resp. difference) module $M$, any tensorial construction on $M$ inherits a unique differential (resp. difference) module structure. In particular, we can see the category of differential (resp. difference) modules as a rigid abelian tensor category (cf. \cite{del}), where the unit object is $K\cdot e$ and $\partial e=0$ (resp. $\Phi e=e$). 
\end{rem}

\begin{defn}
Given an $n$-dimensional differential (resp. difference) module $M$ over $K$, a Picard-Vessiot extension for $M$ (or for (\ref{eqn}), see Remark \ref{rem1}) is a differential (resp. difference) $K$-algebra $R$ of $K$ such that:
\begin{itemize}
\item[a)] the constants in $R$ coincide with the constants in $K$;
\item[b)] the differential module $R\otimes_K M\simeq R^n$ is generated over $R$ by horizontal elements; and,
\item[c)] if $\{f_j=f^i_j\otimes e_i\}_{1\le j\le n}$ are horizontal elements generating $R\otimes_K M$ over $R$, then $R=K[f^i_j,\frac{1}{\det(f^i_j)}]_{1\le i,j\le n}$.
\end{itemize}
\end{defn}

A Picard-Vessiot extension for equation $(\ref{eqn})$ can be obtained as follows. Take the ring of polynomials with coefficients in $K$ and $n\times n$ variables, $K[X^i_j]_{1\le i,j\le n}$, where $n$ is the dimension of $M$. We turn it into a differential or a difference $K$-algebra by setting
\[
v(X^i_j)=a^i_kX^k_j
\]
and we extend it to $K[X^i_j,\frac{1}{W}]_{1\le i,j\le n}$ where $W=\det(X^i_j)_{1\le i,j\le n}$, by putting $v\left(\frac{a}{b}\right)=\frac{a'b-ab'}{b^2}$  if $v$ is a derivation and $v\left(\frac{a}{b}\right)=\frac{v(a)}{v(b)}$ when it is an automorphism, for any $a,b\in K[X^i_j]_{1\le i,j\le n}$. Finally, to obtain a Picard-Vessiot it suffices to form the quotient of $K[X^i_j,\frac{1}{W}]_{1\le i,j\le n}$ by a maximal $v$-stable ideal $I$. For a rigorous exposition of the details involved one could, for example, read ~\cite[Section 1.3]{SvdP} and ~\cite[Section 1.1]{SvdP0}.

Each Picard-Vessiot extension has an associated Galois group $G$. This group is an algebraic group over the field of constants $C$. For our purposes it can be described as follows: the group $GL_n(C)$ acts from the right on $K[X^i_j,\frac{1}{W}]_{1\le i,j\le n}$ by v-commuting $K$-automorphisms by setting
\[
(g^a_b)_{1\le a,b\le n}: X^i_j\longmapsto X^i_kg^k_j
\]
The Galois group $G$ is then the stabilizer of the maximal $v$-stable ideal $I$. Again, full details are in \cite[Section 1.4]{SvdP} and \cite[Section 1.2]{SvdP0}.

The remainder of this article is devoted to proving the following theorem:
\begin{thm}\label{thm}
In the general context, if $G$ is reductive and its determinant group is finite then $I$ is the radical of the ideal generated by the $G$-invariants it contains.
\end{thm}
which is an extension of a theorem by E. Compoint:

\begin{thm}[Compoint]\label{thmcom}
In the differential context with $X=\mathbb{P}^1(C)$, if $G$ is reductive and unimodular then the ideal $I$ is generated by its $G$-invariants.
\end{thm}

The most important aspect of this extension of Compoint's theorem, cf. Corollary \ref{cor} below, is that we can use it to obtain an explicit, finite, simple and useful set of generators of $I$.

\section{The proof}

We start by fixing some notation. Given a $g=(g^i_j)_{1\le i,j\le n}\in GL_n(C)$, the image by the above-mentioned right action on $P\in K[X^i_j,\frac{1}{W}]_{1\le i,j\le n}$ will be denoted by $P^g$. Similarly, there is left action:
\[
(g^a_b)_{1\le a,b\le n}: X^i_j\longmapsto \sum_k g^k_iX^k_j
\]
and the image of $P$ will be denoted by ${}^gP$. We also extend left and right action and their corresponding notation to the field $K(X^i_j)_{1\le i,j\le n}$

\begin{thm}[First main theorem of Invariant Theory]
Let $P\in C[X^i_j]_{1\le i,j\le n}$ and assume that there is a homomorphism $\chi: GL_n(C)\rightarrow C^*$ such that for every $g\in GL_n(C)$, $P^g=\chi (g)P$. Then $P=\lambda \det(X^i_j)_{1\le i,j\le n}^N$ for some $N\in\mathbb{N}$ and some $\lambda\in C$; furthermore $\chi (g)=\det (g)^N$.
\end{thm}

\noindent\emph{Proof:} We may start by assuming that $P$ is homogeneous. As $\chi$ is a homomorphism of $GL_n(C)$ into $C^*$, $\chi(g)=\det(g)^{N_0}$ for some $N_0\in\mathbb{Z}$. Hence $P$ is homogeneous and invariant under $SL_n(C)$, so $P=\lambda \det(X^i_j)_{1\le i,j\le n}^N$. The proof of these statements may be found in~\cite[Paragraph 15]{fulton}. We explain now why the hypothesis forces the polynomial to be homogeneous. Indeed, we just proved that the theorem can be applied to each homogeneous component of $P$; so that under the action of $g\in GL_n(C)$, each component is multiplied by $\det(g)^{N_0}$ for some $N_0$ which increases as the degree gets greater, whence $P$ has a unique component.\hfill$\bigstar$

\begin{cor}
Let $r\in K(X^i_j)_{1\le i,j\le n}$ and assume that there is a homomorphism $\chi: GL_n(C)\rightarrow K^*$ such that for every $g\in GL_n(C)$, $r^g=\chi (g)r$. Then $r=\lambda \det(X^i_j)_{1\le i,j\le n}^N$ for some $N\in\mathbb{Z}$ and some $\lambda\in K$; furthermore $\chi (g)=\det (g)^N$.
\end{cor}

\noindent\emph{Proof:} Let $P,Q\in K[X^i_j]_{1\le i,j\le n}$ be such that $r=P/Q$ where $P$ and $Q$ have no common factors. From the equality $r^g=\chi (g)r$, we obtain the relation $PQ^g=\chi(g)P^gQ$, which says that $P$ divides $P^gQ$. As $K[X^i_j]_{1\le i,j\le n}$ is a unique factorization domain, and $P$ and $Q$ are co-prime, we conclude that $P$ divides $P^g$. Symmetrically, $P^g$ divides $P$ and so  $P^g=\chi_P(g)P$ for some $\chi_P(g)\in K$. Similarly $Q^g=\chi_Q(g)Q$ for some $\chi_Q(g)\in K$. So in order to complete the proof, it suffices to prove the result for $r\in K[X^i_j]_{1\le i,j\le n}$.

Let $\Upsilon$ be a discrete valuation of $K$ such that $\Upsilon(c)=0$ for every coefficient $c$ of $r$. Whence $r$ has coefficients in the valuation ring $K_\Upsilon\subseteq K$ of $\Upsilon$. Let $p$ be the maximal ideal of $K_\Upsilon$, and denote by $r_p$, $r^g_p$ the image of $r$ and $r^g$ in $K_\Upsilon/p\ [X^i_j]$ and $\chi_p (g)$ the image of $\chi(g)$ in $K_\Upsilon/p\simeq C$. We now have $r_p\in C[X^i_j]_{1\le i,j\le n}$, and for every $g\in GL_n(C)$ there is a $\chi_p (g)\in C^*$ such that $r_p^g=\chi_p (g)r_p$. The previous theorem implies that $r_p=\lambda_p \det(X^i_j)_{1\le i,j\le n}^N$ for some $N\in\mathbb{N}$ and some $\lambda_p\in C$. Because every coefficient of $r$ is a unit in $K_\Upsilon$, there is a correspondence between the non-zero coefficients of $r_p$ and the non-zero coefficients of $r$. That is, the monomials of $r$ are the monomials of $\det(X^i_j)$ multiplied by elements of $K_\Upsilon$. Moreover, permuting the monomials of $r$ using permutation matrices of $GL_n$ we conclude that all the coefficients on the monomials coincide. Thus, $r=\lambda\det(X^i_j)$ for some $\lambda\in K$.\hfill$\bigstar$

\begin{prop}
Let $P_1,\ldots, P_r\in K[X^i_j,\frac{1}{W}]_{1\le i,j\le n}$ be linearly independent over $K$. Moreover suppose that they generate an $r$-dimensional $K$-vector space stable under the right $GL_n(C)$-action. Then there exist $g_1,\ldots, g_r\in GL_n(C)$ such that
\[
\det({}^{g_i}P_j)_{1\le i,j\le r}=\lambda\det(X^i_j)^N_{1\le i,j\le n}
\]
for some $\lambda\in C^*$ and $N\in\mathbb{Z}$. Furthermore, one can choose $g_1=e$, the identity element.
\end{prop}

\noindent\emph{Proof:} Consider the vector space $\left(K(X^i_j)_{1\le i,j\le n}\right)^r$. Denote by $s$ the $K(X^i_j)_{1\le i,j\le n}$-rank of the set of vectors
\[
\left\{({}^gP_1,\ldots,{}^gP_r)\in \left(K(X^i_j)_{1\le i,j\le n}\right)^r|\ g\in GL_n(C)\right\},
\]
and choose $g_1,\ldots, g_s$ such that the vectors $({}^{g_i}P_1,\ldots,{}^{g_i}P_r)$, $1\le i\le s$, are linearly independent. Replacing $g_i$ by $g_1^{-1}g_i$ we may assume $g_1=e$.\\
We claim that $s=r$. Assume the contrary, that $s<r$. As $g_1=e$, the determinant of the $(s+1)\times(s+1)$-matrix
\[
\left(\begin{array}{ccc}
P_1         & \cdots & P_{s+1}\\
{}^{g_1}P_1 & \cdots & {}^{g_1}P_{s+1}\\
\vdots      &        & \vdots\\
{}^{g_s}P_1 & \cdots & {}^{g_s}P_{s+1}\\
\end{array}\right)
\]
is zero. Expanding this determinant along the top row we have
\[
P_1\Delta_1-P_2\Delta_2+\ldots+(-1)^{s}P_{s+1}\Delta_{s+1}=0.
\]
By re-indexing the $P_i's$ if necessary, we may assume that the last minor $\Delta_{s+1}=\det({}^{g_i}P_j)_{1\le i,j\le s}$ does not vanish, and we obtain
\[
\frac{\Delta_1}{\Delta_{s+1}}P_1-\frac{\Delta_2}{\Delta_{s+1}}P_2+\ldots+(-1)^{s-1}\frac{\Delta_s}{\Delta_{s+1}}P_s=-(-1)^{s}P_{s+1}
\]

Since the left $GL_n(C)$-action acts on the rows of our $(s+1)\times(s+1)$-matrix, $\frac{\Delta_i}{\Delta_{s+1}}$ is invariant, implying that these quotients of minors are actually in $K$. This contradicts our hypothesis on the linear independence of $P_1,\ldots, P_r$. So $s=r$, and $\det({}^{g_i}P_j)_{1\le i,j\le r}\ne 0$.

The right and the left $GL_n(C)$-actions commutes, and since the $K$-vector space spanned by $P_1,\ldots,P_r$ is invariant under the right $GL_n(C)$-action, we conclude that the right  $GL_n(C)$-action acts on the columns of $\det({}^{g_i}P_j)_{1\le i,j\le r}$, so that for every $g\in GL_n(C)$ there is a $\chi (g)\in K^*$ such that
\[
(\det({}^{g_i}P_j)_{1\le i,j\le r})^g=\chi (g)\det({}^{g_i}P_j)_{1\le i,j\le r}
\]
The corollary above implies now this proposition.\hfill$\bigstar$

\bigskip

Before we prove our theorem, we need the following interpretation of the Galois Correspondence.

\begin{lem}
If $P\in K[X^i_j,\frac{1}{W}]_{1\le i,j\le n}$ is $G$-invariant, then there is an $f\in K$ such that $P(x)-f\in I$.
\end{lem}

\noindent\emph{Proof}: Under the canonical map $K[X^i_j,\frac{1}{W}]_{1\le i,j\le n}\mapsto K[X^i_j,\frac{1}{W}]_{1\le i,j\le n}/I$, the image of $P$ is invariant under $G$, so by Galois correspondence $P$ maps to an $f$ in to the ground field $C$. Thus, the image of $P-f$ under the canonical map is zero.\hfill$\bigstar$

\bigskip

\noindent\emph{Proof of Theorem \ref{thm}:} Let $Q\in K[X^i_j]_{1\le i,j\le n}$. As the $GL_n(C)$-action does not modify the degree, the degree of the elements in the orbit of $Q$ is bounded; and so, they spans a finite dimensional $C$-vector space. The same holds if $Q \in K[X^i_j,\frac{1}{W}]_{1\le i,j\le n}$, for $W^g=\det(g)W$.

Now take $Q\in I$. We will prove that a power of $Q$ lies in the ideal generated by the $G$-invariants in $I$. Denote by $Q_0^{GL_n(C)}$ the $K$-vector space spanned by $\{Q^g\}_{g\in GL_n(C)}$ and by $r_0$ its dimension; similarly, denote by $Q_0^G$ the space spanned by $\{Q^g\}_{g\in G}$ and by $t_0$ its dimension. Note that $Q_0^G\subseteq I$ as $I$ is $G$-stable. Let $Q_1,\ldots, Q_{r_0}$ be a basis of $Q_0^{GL_n(C)}$ such that $Q_1=Q$ and $Q_1,\ldots,Q_{t_0}$ is a basis of $Q_0^G$. Additionally, as $G$ is reductive we can assume that $Q_{t_0+1},\ldots, Q_{r_0}$ span a $G$-stable space. Indeed, this last condition can be satisfying as every linear representation of a reductive group over $C$ is completely reducible.

Let $N_0\gg 0$ be such that if $g\in G$ then $\det(g^{N_0})=1$. Denote by $Q^{GL_n(C)}$ the space spanned by the $N_0$-fold products of the $Q_1,\ldots, Q_{r_0}$ and by $r$ its dimension; similarly denote by $Q^G$ the sub-space spanned by the $N_0$-fold products of the $Q_1,\ldots, Q_{t_0}$ and by $t$ its dimension. Note that $Q^G\subseteq I$. Let $P_1,\ldots, P_r$ be a basis of $Q^{GL_n(C)}$ such that $P_1=Q^{N_0}$, $P_1,\ldots,P_r$ is a basis of $Q^G$ and the space spanned by $P_{t+1},\ldots, P_r$ is $G$-stable. The previous proposition implies that there exist $g_1=e,g_2,\ldots,g_r\in GL_n(C)$ such that
\[
\det({}^{g_i}P_j)_{1\le i,j\le r}=\lambda\det(X^i_j)^N_{1\le i,j\le n}
\]
for some $\lambda\in K^*$ and $N\in\mathbb{Z}$. Notice that the action of $GL_n(C)$ on $Q^{GL_n(C)}$ is the $N_0$-th symmetric power of the action on $Q_0^{GL_n(C)}$, thus $N$ is a multiple of $N_0$ and $\det(g)^N=1$.

As $\det(X^i_j)_{1\le i,j\le n}=W$ is a unit in $K[X^i_j,\frac{1}{W}]_{1\le i,j\le n}$, then $\det(X^i_j)_{1\le i,j\le n}\not\in I$. Now, expanding the determinant $\det({}^{g_i}P_j)_{1\le i,j\le r}$ through the first row, as $P_1,\ldots,P_t\subseteq I$, we can see that $t<r$ (or else $t=r$ and we would have $W^N\in I$).\\
We compute $\det({}^{g_i}P_j)_{1\le i,j\le r}$ by expanding through its first $t$ columns: given a collection of $t$ indices $J=\{j_1,\ldots,j_t\}$, with $1\le j_1<j_2<\ldots<j_t\le r$, we denote by $P_J$ the determinant
\[
P_J=\det({}^{g_{j_k}}P_{j})_{1\le j,k\le t}
\]
and by $Q_J$ its co-determinant so that
\[
\det({}^{g_i}P_j)_{1\le i,j\le r}=\sum_JP_JQ_J
\]
Note that $P_J\in I$ if $1\in J$. But, as $\det({}^{g_i}P_j)_{1\le i,j\le r}\not\in I$, there is a collection of indices $k=\{k_1,\ldots,k_t\}$ such that $1\not\in k$. Because $P^G$ is $G$-stable, and the left and right actions commute, for every $g\in G$ there is a $\chi(g)\in K$ such that $P_J^g=\chi(g)P_J$ for every $J$. Similarly, because the space spanned by $P_{t+1},\ldots, P_r$ is $G$-stable too, for every $g\in G$ there is a $\chi'(g)\in K$ such that $Q_J^g=\chi'(g)Q_J$ for every $J$. Thus, from the expansion of $\det({}^{g_i}P_j)_{1\le i,j\le r}$ we read that $\chi(g)\chi'(g)=\det(g)^N=1$. In particular, for any two collection of $t$ indices $J$ and $L$, $P_JQ_L$ is $G$-invariant.

The previous lemma implies that there is an $f\in K$ such that $P_kQ_k-f\in I$. Since $P_kQ_k\not\in I$ we know that $f\ne 0$.

Set $k_0=1$, $P_0=P_1=Q^N$. Denote by $k(s)$ the sequence $\{k_0,k_1,\ldots,\hat{k_s},\ldots,k_t\}$ where the hat means that we removed the $s$-th element. In particular $k(0)=k$. Also, the determinant of the $(t+1)\times(t+1)$ matrix
\[
\left(\begin{array}{cccc}
{}^{g_{k_0}}P_0 & {}^{g_{k_0}}P_1 & \cdots & {}^{g_{k_0}}P_{t}\\
{}^{g_{k_1}}P_0 & {}^{g_{k_1}}P_1 & \cdots & {}^{g_{k_1}}P_{t}\\
\vdots          & \vdots          &        & \vdots\\
{}^{g_{k_t}}P_0 & {}^{g_{k_t}}P_1 & \cdots & {}^{g_{k_t}}P_{t}\\
\end{array}\right)
\]
is zero as the two first columns are the same. Expanding along the first column we have:
\[
0=Q^{N_0}P_k-{}^{g_{k_1}}Q^{N_0}P_{k(1)}+\ldots+(-1)^t\ {}^{g_{k_t}}Q^{N_0}P_{k(t)}
\]
Multiplying first by $Q_k$ on both sides and then subtracting $fQ^{N_0}$ we obtain:
\[
-fQ^{N_0}=Q^{N_0}(P_kQ_k-f)-{}^{g_{k_1}}Q^{N_0}P_{k(1)}Q_k+\ldots+(-1)^t\ {}^{g_{k_t}}Q^{N_0}P_{k(t)}Q_k
\]
Finally, note that if $s>0$ then $1\in k(s)$ and so $P_{k(s)}\in I$. Whence we have an expression of $Q^{N_0}$ as a linear combination of $G$-invariants in $I$.\hfill$\bigstar$

\begin{cor}\label{cor}
Assuming the context of Theorem \ref{thm}, suppose $P_1,\ldots,P_r$ generate the $C$-algebra of $G$-invariants in $C[X^i_j]_{1\le i,j\le n}$, and $f_1,\ldots,f_r\in K$ are such that $P_1-f_1,\ldots,P_r-f_r\in I$. Then $I$ is the radical of the ideal generated by $\{P_i-f_i\}_{i\in\{1,\ldots,r\}}$.
\end{cor}

\noindent\emph{Proof:} Because $G$ is reductive, the $C$-algebra of $G$-invariants in $C[X^i_j]_{1\le i,j\le n}$ is finitely generated. Now, if $P\in C[X^i_j,\frac{1}{W}]_{1\le i,j\le n}$ then, for some $N\gg 0$, $W^NP\in K[X^i_j]_{1\le i,j\le n}$. Additionally, if $P$ is $G$-invariant, and $N$ is such that $W^N$ is $G$-invariant, $W^NP$ can be written as a $K$-linear combination of $G$-invariants in $C[X^i_j]_{1\le i,j\le n}$.

We now prove by induction on the degree of $Q\in K[X^i_j]$ that if $Q\in K[X^i_j]^G\cap I$ then $Q\in\langle P_i-f_i\rangle_{i\in\{1,\ldots,r\}}$, the case of degree zero being trivial. Let $Q=Q_{\deg Q}+Q_{\deg Q-1}+\ldots+Q_1+Q_0$ be the homogeneous decomposition of $Q$. Since the $G$-action preserves the degree and $Q$ is $G$-invariant then so is $Q_{\deg Q}$. Let $q\in K$ be such that $Q_{\deg Q}-q\in I$. Then $Q-(Q_{\deg Q}-q)\in I$ and $\deg (Q-(Q_{\deg Q}-q))< \deg Q$, so by induction hypothesis $Q-(Q_{\deg Q}-q)\in\langle P_i-f_i\rangle_{i\in\{1,\ldots,r\}}$. It remains to prove that $(Q_{\deg Q}-q)\in \langle P_i-f_i\rangle_{i\in\{1,\ldots,r\}}$. Note that if $Q_1-q_1, Q_2-q_2\in I$ then $(Q_1+Q_2)-(q_1+q_2)\in I$ and $Q_1Q_2-q_1q_2=(Q_1-q_1)(Q_2-q_2)+q_2(Q_1-q_1)+q_1(Q_2-q_2)\in I$. Thus if $\Phi$ is a polynomial over $K$ with $r$ variables such that $Q_{\deg Q}=\Phi(P_1,\ldots, P_r)$ then $q=\Phi(f_1,\ldots,f_n)$ and the statement of the last observation implies $\Phi(P_1,\ldots,P_r)-\Phi(f_1,\ldots,f_r)\in I$.\hfill$\bigstar$

\end{document}